\def\Aa{{\E{A}}}
\def\BB{{\E{B}}}
\def\DD{{\E{D}}}
\def\al{\alpha}
\def\be{\beta}
\def\de{\delta}
\def\la{\lambda}
\def\De{\Delta}
\def\vep{\varepsilon}
\def\oli{\overline}
\def\sbs{\subset}
\def\sbst{\subseteq}
\def\no={\,{\,|\!\!\!\!\!=\,\,}}
\def\Raw{\Rightarrow}
\def\rank{\text{rank}}
\def\wt{\widetilde}
\def\wh{\widehat}
\def\ssm{\smallsetminus}

\documentclass[12pt,final]{amsart}
\usepackage{amssymb}
\usepackage{lpretex}

\DocVersion[Close to final, T.E.]{5}
\DocVersion[Fixed several things mainly the indexing of cohomology, T.E.]{6}
\DocVersion[Tried to correct thm 4.8 and 4.9, A.B.; Nov. 26]{7}
\DocVersion[Some further mod's of 4.8 and 4.9, T.E.; Nov. 28]{8}

\def\PRcanbb#1{\ensuremath{{\mathbb #1}}}

\theoremstyle{plain}
\newtheorem{Thm}{Theorem}[section]
\newtheorem{Lem}[Thm]{Lemma}
\newtheorem{Prop}[Thm]{Proposition}

\numberwithin{equation}{section}

\def\bpf{\begin{proof}}
\def\epf{\end{proof}}

\DeclareMathAlphabet{\E}{U}{eus}{m}{n}{}

\theoremstyle{remark}

\def\DHitemrefadd#1{(#1)}
\begin{document}

\title[Subspace arrangements over finite fields]
{Subspace arrangements over finite fields: Cohomological and 
enumerative aspects}

\author[Anders Bj\"orner]{Anders Bj\"orner}
\address{\hskip-\parindent Anders Bj\"orner\\
 Department of Mathematics\\
 Kungl.~Tekniska H\"ogskolan\\
 S-100 44  Stockholm\\
 Sweden} 
\email{bjorner@math.kth.se}
\author[Torsten Ekedahl]{Torsten Ekedahl}
\thanks{Bj\"orner was partially supported by the Institute for
Advanced Study (Princeton, NJ) and the 
Mathematical Sciences Research Institute (Berkeley, CA).
Ekedahl was 
supported by the G\"oran Gustafsson Foundation for Research in
Natural Sciences and Medicine}
\address{\hskip-\parindent Torsten Ekedahl\\Department of Mathematics\\
 Stockholm University\\
 S-106 91  Stockholm\\
Sweden}
\email{teke@matematik.su.se}

\begin{abstract}
The enumeration of points on (or off) the union of some linear or affine subspaces
over a finite field is dealt with in combinatorics via the characteristic
polynomial and in algebraic geometry via the zeta function. We discuss the
basic relations between these two points of view. Counting points is also
related to the $\ell$-adic cohomology of the arrangement (as a variety).
We describe the eigenvalues of the Frobenius map acting on this cohomology,
which corresponds to a finer decomposition of the zeta function. The
$\ell$-adic cohomology groups and their decomposition into eigenspaces are
shown to be fully determined by combinatorial data. Finally, it is shown
that the zeta function is determined by the topology of the corresponding
complex variety in some important cases.
\end{abstract}
\maketitle

\begin{section}{Introduction}

This paper is concerned with subspace arrangements over general (in
particular finite) fields, and with their enumerative and cohomological
invariants. In this introduction we will summarise the main
results. We begin with a review of some algebraic geometry.

Let $V$ be a $d$-dimensional projective variety over the field $\F_q$ of
order $q=p^\al$. (We will not assume that a variety is irreducible.)
For each extension field $\F_{q^s}$, $s\ge1$, let $N_s$ be 
the number of points on $V$ over $\F_{q^s}$. The zeta function of $V$ is
the formal power series
$$
Z(V; t)=\exp\left(\sum_{s\ge1} N_s\,\frac{t^s}{s}\right).
$$

Let $\ell$ be a prime number, $\ell\no=p$, and $\Q_\ell$ the field of
$\ell$-adic numbers. Let $H^i (V, \Q_\ell)$ be the $\ell$-adic \'etale
cohomology groups of Grothendieck \cite{D4,G4,G5}. These are
finite-dimensional vector spaces over $\Q_\ell$, a field of
characteristic zero. Furthermore, $H^i(V, \Q_\ell)=0$ unless $0\le
i\le 2d$.

The Frobenius map \pil FVV, defined by
$(x_1,\dots,x_n)\mapsto\left(x^q_1,\dots,x^q_n\right)$, induces a linear 
action
$F:H^i(V,\Q_\ell)\to H^i(V,\Q_\ell)$ for each $0\le i\le 2d$. Let
$\al_{ij}$ be the eigenvalues of this map. The Grothendieck-Lefschetz
fixed point formula \cite{G5} 
implies that
$N_s=\sum_{i,j}(-1)^i\al^s_{i,j}$, which is equivalent to the following
decomposition of the zeta function:
$$
Z(V;t)=\frac{P_1(t) P_3(t)\dots P_{2d-1}(t)}{P_0(t) P_2(t)\dots
P_{2d}(t)}
\qquad,\qquad\text{ where\,} P_i(t)=\prod_j(1-\al_{ij}t).
$$
Note that $P_i(t)\in\Q_\ell[t]$, and this polynomial might, seen a
priori, depend on $\ell$. However, this is not the case and the
coefficients of $P_i(t)$ are in fact algebraic integers. Deligne
\cite{D1} showed that $|\al_{ij}|=q^{i/2}$ if $V$ is smooth.

In this paper we determine the polynomials $P_i(t)$ for the case when
$V$ is a union of linear subspaces. It  turns out that $P_i(t)$ is
then determined by combinatorial data in the following way.

Let $\Aa=\{K_1,\dots,K_t\}$ be an arrangement of linear subspaces in projective
$(n-1)$-space over $\F_q$, and let $V_\Aa$ denote their union (a singular
projective variety). Let $L_\Aa$ denote the partially ordered set of all
nonempty intersections $K_{i_1}\cap\dots\cap K_{i_m}$, $1\le i_1<\dots<i_m\le
t$, ordered by reverse inclusion. Let $\be^{\ge j}_i$ denote the rational
$i$-th Betti number of the (simplicial) order complex homology of the subposet $\{x\in
L_\Aa\mid j\le\dim(x)<n-1\}$.

\begin{Thm}\label{J}
For the union of a $d$-dimensional projective subspace arrangement
$\Aa$ over $\F_q$ we have
$$
P_i(t)=\prod^d_{j=0}\left(1-q^jt\right)^{\be^{\ge j}_{i-2j}}.
$$
\end{Thm}

The formula can be seen as a sharpening or Frobenius-equivariant
version of a recent formula for $H^i(V_\Aa,\Q_\ell)$ due to Yan
\cite{Ya}. Theorem \ref{J} will be proved with a cohomological
argument that is valid for arbitrary fields. In particular, a unified
setting will be given for earlier results such as Yan's theorem and
the complex version of Ziegler and \v{Z}ivaljevic \cite{ZZ}.

One can also consider $\ell$-adic cohomology and eigenvalues of
Frobenius on the {\em complement} (rather than the union) of a
projective arrangement. Or one could consider these questions for
arrangements of affine subspaces of affine $n$-space over $\F_q$. Of
the four possible combinations, affine/projective and
union/complement, we have chosen here to concentrate on unions in the
projective case and complements in the affine case. Formulas covering all
cases, and also punctured affine arrangements, are given in the paper.

We will now state the result in the affine case. Let $\Aa$ be an
arrangement in $\A^n=\F^n_q$ with union $V_\Aa$ and complement
$M_\Aa=\F^n_q\setminus V_\Aa$. Let $\al_{ij}$ be the eigenvalues of the induced
Frobenius map on $\ell$-adic cohomology with compact supports
$F:H^i_c(M_\Aa, \Q_\ell)\to H^i_c(M_\Aa, \Q_\ell)$. Let
$P_i(t)=\prod_j(1-\al_{ij}t)$ as before. Let $\be^{\oplus
j}_i=\sum\be_i(\hat0,x)$, where $\be_i(\hat0,x)$ denotes the
$i$-th rational Betti number of the order complex homology of
the open interval $(\hat0,x)$ in $L_\Aa$, and the sum is over all
$x\in L_\Aa$ such that $\dim(x)=j$.

\begin{Thm}\label{K}
For the complement of a $d$-dimensional affine subspace arrangement
$\Aa$ over $\F_q$ we have
$$
P_i(t)=\prod^d_{j=0}\left(1-q^jt\right)^{\be^{\oplus j}_{i-2j-2}}.
$$
\end{Thm}
A corresponding decomposition of $H^i_c (M_\Aa, \Q_\ell)$ without the
eigenvalue information was given by Yan \cite{Ya}, and in the real
case earlier by Goresky and MacPherson \cite{GM}. For the special case
of hyperplane arrangements Theorem \ref{K} specialises to say that
Frobenius acts on $H^i_c(M_\Aa, \Q_\ell)$, for $n\le i\le 2n$, 
with only one eigenvalue, namely $q^{i-n}$, and for all other $i$
these cohomology groups vanish; a result earlier obtained by Lehrer \cite{Le} 
(see also Kim \cite{Ki}). 

The paper is organised as follows: Some definitions and facts about
subspace arrangements and the combinatorics of intersection
semilattices are reviewed in Section 2. In Section 3 we discuss
counting points on (or off) subspace arrangements over finite fields.
Some formulas relating zeta functions to the characteristic polynomials of
intersection semilattices are given along with some related facts. The proofs of the main
cohomological results outlined above are given in Section 4. Specifically,
Theorems \ref{J} and \ref{K} appear as part of Theorems \ref{coharr} and 
\ref{cohcomp}. For
arrangements $\Aa$ defined by forms with integer coefficients we can
consider both the arrangement $\Aa_q$ over $\F_q$ ($q=p^\al$),
obtained by reduction modulo $p$, and the complex arrangement
$\Aa_\C$. In the final section we discuss connections between the zeta
function of $\Aa_q$ and the Betti numbers of $\Aa_\C$, showing that
in some important cases they mutually determine each other. 

Valuable conversations with P. Deligne and K.S. Sarkaria are 
gratefully acknowledged.
\end{section}

\begin{section}{Preliminaries}

Let $\F$ be a field. By an {\em affine subspace arrangement} we mean a
finite collection of affine subspaces in $\F^n$. Similarly, by a {\em
projective subspace arrangement} we mean a finite collection of linear
subspaces in projective $(n-1)$-space $\F\P^{n-1}$. An arrangement of either
kind is {\em essential} if the intersection of all subspaces is
empty. An affine arrangement is called {\em central} if all subspaces
contain the origin. There is an obvious one-to-one correspondence
between central arrangements in $\F^n$ and projective arrangements in
$\F\P^{n-1}$. An arrangement (of either kind) is {\em $d$-dimensional}
if $d$ is the maximal dimension of one of its subspaces.

For an arrangement $\Aa=\{K_1,\dots,K_t\}$ we denote by $V_\Aa$ the
{\em variety} $K_1\cup\dots\cup K_t$. Also, we let
$L_\Aa=\{K_{i_1}\cap\dots\cap K_{i_m}\no=\varnothing\}$ be the {\em
intersection semilattice}, the family of all nonempty intersections of
subarrangements ordered by reverse inclusion. The semilattice $L_\Aa$
has a least element $\hat0$ (equal to the whole space $\F^n$ or
$\F\P^{n-1}$, as the case may be), and if $\Aa$ is central (or equivalent
to a central arrangement) then there is also a greatest element
$\hat1$ (equal to $K_1\cap\dots\cap K_t$).

We refer to \cite{Bj} for a general introduction to the theory of
subspace arrangements.

Let $P$ be a poset (short for ``partially ordered set'')
and $x,y\in P, x<y$. Then $[x,y]=\{z\in P\mid x\le
x\le y\}$ and $(x,y)=\{z\in P\mid x<z<y\}$ are the {\em closed} and
{\em open intervals}.  Also, for $p \in P$ let $P_{<p}$ resp.~$P_{\le p}$ 
be the subset consisting of those elements of $P$ which are
less than $p$ (resp.~less than or equal to $p$).
We assume familiarity with the {\em M\"obius
function} $\mu(x,y)$ of $P$, see \cite{S2}.

With a poset $P$ we associate its {\em order complex} $\De(P)$
consisting of all chains $x_0<x_1<\dots<x_k$. This is a simplicial
complex, so we obtain (simplicial) homology groups $H_i(P)=H_i(\De(P);
\Z)$, Betti numbers $\be_i(P)=\rank_\Z H_i(P)$, Euler characteristic
$\chi(P)$, etc.

A tilde will always denote reduced homology (Betti numbers, Euler
characteristic): $\wt{H}_i, \wt\be_i, \wt\chi$. Recall Ph. Hall's
theorem \cite[p. 120]{S2}.

\begin{equation}\label{a}
\mu(x,y)=\wt\chi(\De(x,y)),
\end{equation}
for all $x<y$ in $P$.

We will sometimes use the following arrangements, the ``$k$-equal
arrangements of type A, B and D'' \cite{BL,BSag,BWe}, to provide examples:
\begin{align*}
\Aa_{n,k}
&=
\{x_{i_1}=x_{i_2}=\dots=x_{i_k}\mid 1\le i_1<\dots<i_k\le n\}\\
\DD_{n,k}
&=
\{\vep_1 x_{i_1}=\vep_2 x_{i_2}=\dots=\vep_k x_{i_k}\mid
1\le i_1<\dots<i_k\le n, \vep_i\in\{+1,-1\}\}\\
\BB_{n,k}
&=
\DD_{n,k}\cup\{x_{j_1}=\dots=x_{j_{k-1}}=0
\mid 1\le j_1<\dots<j_{k-1}\le n\}.
\end{align*}
\end{section}

\begin{section}{Characteristic polynomial and zeta function}

In this section we will develop some of the basic facts about counting
points on (or off) arrangements over finite fields. The tool for this
in combinatorics is the characteristic polynomial of the arrangement
and in algebraic geometry the zeta function. We will make the
relationship between these notions explicit and prove a few related
facts.

The following is a version of the combinatorial ``principle of
inclusion-exclusion''.

\begin{Prop}\label{B} Let $\{H_i\}_{i\in I}$ be a family of subsets of
a finite set $E$. Let $L$ be the semilattice of nonempty intersections
of subfamilies, ordered by reverse inclusion. Then
$$
\text{\rm{card}} \left(E-\bigcup_{i\in I} H_i\right)=\sum_{x\in L} \mu\left(\hat0,x\right) 
\text{\rm{card}}(x).
$$
\end{Prop}

\bpf For each $e\in E$ let $x_e$ be the intersection of all $H_i$
containing $e$. If $e\in H_i$ for some $i\in I$, then $x_e\no=\hat0$
and the right hand side will count $e$ altogether $\sum_{\hat0\le x\le
x_e} \mu\left(\hat0,x\right)=0$ times. Otherwise, $x_e=\hat0$ and $e$ will be
counted $\mu\left(\hat0,\hat0\right)=1$ time.
\epf

For an affine arrangement $\Aa$ in $\F^n$ let
\begin{equation}\label{d}
P_\Aa(t)=\sum_{x\in L_\Aa} \mu\left(\hat0,x\right)t^{\dim(x)}
\end{equation}
define the {\em characteristic polynomial} of $\Aa$. Such polynomials
have long been studied for hyperplane arrangements \cite{OT,Za} and
for other graded posets \cite{S2}; they were generalised to subspace
arrangements in \cite{Bj,BL}.

\begin{Prop}\label{C} 
Let $q$ be a prime power, and let $\Aa$ be an affine subspace
arrangement in $\F^n_q$. Then:
$$
P_\Aa(q)=\text{\em{card}\,}\left(\F^n_q\setminus V_\Aa\right).
$$
\end{Prop}

\bpf If $x$ is an affine subspace then $\text{card\,}(x)=q^{\dim(x)}$,
so this is a special case of Proposition \ref{B}.
\epf

\begin{remark}
This result is well-known for hyperplane arrangements, see the
``critical problem'' in \cite{CR} and also \cite[p. 51]{OT}. The
extension to subspace arrangements has independently been found by
Athanasiadis \cite{At}. A similar result (in a somewhat different
setting) appears in Blass and Sagan \cite{BlS}.
\end{remark}
Now, let $\Aa$ be a projective arrangement in $\F\P^{n-1}$ and define
\begin{equation}\label{l}
P^*_\Aa(t)=\sum_{x\in L_\Aa} \mu\left(\hat0,x\right)
\left(1+t+t^2+\dots+t^{\dim(x)}\right).
\end{equation}

\begin{Prop}\label{D} Let $\Aa$ be a projective arrangement in
$\F_q\P^{n-1}$. Then:
$$
P^*_\Aa(q)=\text{\em{card}\,}\left(\F_q\P^{n-1}\setminus V_\Aa\right).
$$
\end{Prop}

\bpf This follows from Proposition \ref{B}, since
$\text{card\,}(x)=1+q+\dots+q^{\dim(x)}$ for projective subspaces
$x$. 
\epf

Note that for projective $\Aa$, if $\oli\Aa$ is the corresponding central arrangement in
$\F^n$, then
\begin{equation}\label{e}
P_{\oli\Aa}(t)=(t-1)P^*_\Aa(t).
\end{equation}
This follows from the preceding propositions, since
every point in the complement
of $\Aa$ corresponds to $q-1$ points in the complement of $\oli\Aa$.

We will call $P^*_\Aa(t)$ the {\em reduced characteristic polynomial}
of a projective arrangement $\Aa$. Its coefficients are (up to sign)
the reduced Euler characteristics of the $j$-{\em truncations} of the
intersection semilattice $L_\Aa$, defined by
$$
L^{\ge j}_\Aa=\{x\in L_\Aa\setminus \{\hat0\}\mid\dim(x)\ge j\}.
$$

\begin{Prop}\label{E} Let $P^*_\Aa(t)=\sum^{n-1}_{j=0} c_j t^j$. Then
$c_j=-\wt\chi\left(L^{\ge j}_\Aa\right)$.
\end{Prop}
\bpf Relations \eqref{a} and \eqref{l} imply that
$$
\wt\chi\left(L^{\ge j}_\Aa\right)=-\!\!\!\sum_{\dim(x)\ge j}
\!\!\mu\left(\hat0,x\right)=-c_j.
$$
\epf
It follows from the definition that the characteristic polynomial of a
$d$-di\-men\-sio\-nal affine arrangement in $\F^n$ has the structure
\begin{equation}\label{j}
P_\Aa(t)=c_0+c_1t+\dots+c_dt^d+t^n,\qquad\text{with\, } c_d<0.
\end{equation}
Similarly, the reduced characteristic polynomial of a $d$-dimensional
projective arrangement in $\F\P^{n-1}$ has the structure
\begin{equation}\label{k}
P^*_\Aa(t)=c_0+c_1t+\dots+c_dt^d+t^{d+1}+\dots+t^{n-1},\qquad\text{with\,
}c_d\le0.
\end{equation}
Hence, there are only $d+1$ essential coefficients $c_0,\dots,c_d$ in
each case.
\begin{example}
The 3-equal arrangement $\Aa=\Aa_{6,3}$ gives
for its 4-dimensional affine version
$$
P_\Aa(t)=-26t^2+45t^3-20t^4+t^6,
$$
and for its 3-dimensional projective version
$$
P^*_\Aa(t)=26t^2-19t^3+t^4+t^5.
$$
\end{example}
Let $\Aa$ be an affine subspace arrangement in $\F^n_q$, for some
prime power $q$. Consider the field extensions
$\F_q\sbs\F_{q^s}\sbs\oli\F_q$, where
$\oli\F_q$ denotes an algebraic closure of
$\F_q$. Then $\Aa$ can also be considered as an arrangement in
$\oli\F^n_q$, defined by equations with coefficients in $\F_q$, and
its intersection with $\F^n_{q^s}$ for all the intermediate fields is
well defined.

For arbitrary subsets $X\sbst\oli\F^n_q$ define the {\em zeta function}
of $X$ as the formal power series
\begin{equation*}
Z(X;t)=Z_q(X;t)=\exp\left(\sum_{s\ge1}\text{card}\left(X\cap
\F^n_{q^s}\right)\frac{t^s}{s}\right).
\end{equation*}

The analogous definition is made in the projective case of the {\em
zeta function} for subsets $X\sbst\F_q\P^{n-1}$:
\begin{equation*}
Z(X;t)=Z_q(X;t)=\exp\left(\sum_{s\ge1}\text{card}
\left(X\cap \F_{q^s}\P^{n-1}\right)\frac{t^s}{s}\right).
\end{equation*}

The following lemma will be useful. It is no doubt well known but
we lack an appropriate reference.

\begin{Lem}\label{F} If $P(t)=\sum^d_{j=0} c_jt^j$, then
$\exp\left(\sum_{s\ge1}P(q^s)\frac{t^s}{s}\right)=\prod^d_{j=0}\left(1-q^j
t\right)^{-c_j}.$
\begin{proof}
Taking the logarithm and the derivative of both sides reduces 
the identity to the summing
of a linear combination of geometric series. The constant lost by taking the
derivative is determined by putting $t$ equal to 0.
\end{proof}
\end{Lem}

There is the following connection between characteristic polynomials
and zeta functions.

\begin{Thm}\label{G} {\rm(i)} Let $\Aa$ be an affine arrangement in
$\F^n_q$, with characteristic polynomial
$P_\Aa(t)=\sum^n_{j=0}c_jt^j$. Then
$$
Z\left(\oli\F^n_q\setminus V_\Aa;t\right)=\prod^n_{j=0}
\left(1-q^jt\right)^{-c_j}.
$$
{\rm(ii)} Let $\Aa$ be a $d$-dimensional projective arrangement in
$\F_q\P^{n-1}$, with reduced characteristic polynomial
$P^*_\Aa(t)=\sum^{n-1}_{j=0} c_jt^j$. Then
$$
Z(V_\Aa;t)=\prod^d_{j=0}\left(1-q^jt\right)^{c_j-1}.
$$
\end{Thm}

\bpf 
Note that the intersection semilattice of $\Aa$ is the same over
all fields $\F_{q^s}$ (cf. the proof of Lemma \ref{A}). Hence, part
(i) is an immediate consequence of Proposition \ref{C} and Lemma
\ref{F}.

Similarly, for part (ii) Proposition \ref{D} gives
$$
\text{card}\left(V_\Aa\cap\F_{q^s}\P^{n-1}\right)=\frac{q^{sn}-1}{q^s-1}
-P^*_\Aa\left(q^s\right)=\sum^{n-1}_{j=0}(1-c_j)q^{sj}.
$$
Now, by \eqref{k} we have that $c_j=1$ for $j=d+1,\dots,n-1$. Hence
the result follows from Lemma \ref{F}.
\epf

The zeta function for the union of a $d$-dimensional projective
arrangement can because of Proposition \ref{E} also be stated in this
form:
\begin{equation}\label{i}
Z(V_\Aa;t)=\prod^d_{j=0}\left(1-q^jt\right)^{-\chi\left(L^{\ge j}_\Aa\right)}.
\end{equation}

It is known from the work of Dwork \cite{Dw} that the zeta function of {\em any}
affine or projective variety over $\F_q$ is a rational power series with
coefficients in $\N$. In the theory of formal power series there is a slightly
stronger concept. A series is called $\N$-{\em rational} if it can be produced
from $\N$-polynomials by finitely many times applying the operations of
addition, multiplication and quasi-inverse (this is defined as $f+f^2+f^3+\dots$
for series $f$ such that $f(0)=0$). The point is that subtraction is never
allowed. See Eilenberg \cite{Ei} and Reutenauer \cite{Re} for more about this
concept. Reutenauer points out that it follows from the work of Deligne
\cite{D1} that the zeta function of a smooth variety is in fact
$\N$-rational. The same can be said about the zeta functions considered here,
indeed Reutenauer's observation is true without the assumption of smoothness.

\begin{Thm}\label{I} The zeta function of any algebraic variety $X$ defined 
over the finite field $\F_q$ is $\N$-rational.
\end{Thm}
\bpf
A theorem of Soittola \cite{So} says that a rational series with
coefficients in $\N$ is $\N$-rational if it has a real pole $\al>0$
such that $|\be|>\al$ for all other poles $\be$. Assume first that $X$ is
\Definition{absolutely irreducible}, i.e., is irreducible over an algebraic
closure of $\F_q$. Then the estimate by Lang and Weil \cite[Theorem 1]{LW}
gives that the number $N_s$ of points of $X$ over an extension field of
cardinality $q^s$ fulfills a uniform estimate $|N_s-q^{ds}|\le
O(q^{s(d-1/2)})$, where $d$ is the dimension of $X$. This gives that
$Z(X;t)-1/(1-q^dt)$ is holomorphic in the disc $|t|<q^{1/2-d}$, and hence all
the poles of $Z(X;t)$ except $q^{-d}$  have absolute value larger than
$q^{-d}$.

In the more general case when $X$ is irreducible but not necessarily absolutely
irreducible, the functions on $X$ algebraic over $\F_q$ form a finite
extension field $\F_{q^n}$. As then $X$ has points over a field $\F_{q^m}$
only if $n|m$ we see that $Z(X;t)=Z(X';t^n)$, where $X'$ is $X$ considered as
a variety over $\F_{q^n}$. Since the substitution $t \mapsto t^n$ preserves
$\N$-rationality we reach the conclusion also in this case.

In the most general case we note that the zeta function is multiplicative over
disjoint unions of subvarieties and that any variety is the disjoint union of
finitely many irreducible varieties.
\epf
\begin{remark}
Instead of using the result of Lang-Weil we could have made an appeal to
Deligne's theorem. However, the former predates the latter and is considerably
more elementary.
\end{remark}

It is a well known fact for affine and projective hyperplane
arrangements, originally due to G.-C.~Rota, that the coefficients of
the corresponding characteristic polynomial ``alternate in sign'',
meaning that $c_j c_{j+1}\le0$ for all $j$,
cf. \cite[p. 126]{S2}. This has via Theorem \ref{G} consequences for
the structure of the zeta function. We will now show that the same is
true for certain subspace arrangements.

Let us say that an intersection semilattice $L_\Aa$ is {\em
hereditary} if whenever $x\in L_\Aa\setminus \{\hat0\}$ and $\dim(x)>0$ there
exists $y\in L_\Aa$ such that $y>x$ and $\dim(y)=\dim(x)-1$. We will
say that $L_\Aa$ is {\em mod-$2$-pure} if either all  maximal chains
are of even length or all maximal chains are of odd length. We will
consider the concept of CL-shellability known, see \cite{BW1,BW2}
for all details.

\begin{Thm}\label{H} Suppose $L_\Aa$ is hereditary, mod-$2$-pure and
CL-shellable.
\begin{enumerate}
\item{If $\Aa$ is affine and $d$-dimensional with characteristic
polynomial $P_\Aa(t)=\sum^d_{j=0}$\break
$ c_jt^j+t^n$, then $(-1)^{d-j}c_j\le0$
for all $0\le j\le d$.}
\item{If $\Aa$ is projective and $d$-dimensional with reduced
characteristic polynomial
$P^*_\Aa (t) = \sum^d_{j=0} c_j t^j + \sum^{n-1}_{j=d+1} t^j$, then
$(-1)^{d-j}c_j\le0$ for all $0\le j\le d$.}
\end{enumerate}
\end{Thm}
\begin{proof}
Suppose that all maximal chains in $L_\Aa$ are of even
length. The case of odd length is handled similarly.
\begin{enumerate}
\item{Let $x\no=\hat0$. The interval $[\hat0,x]$ is CL-shellable
\cite[Lemma 5.6]{BW1}, and being hereditary implies that the length of
any maximal chain in $[\hat0,x]$ has the same parity as
$\dim(x)$. Hence, by \cite[Proposition 5.7]{BW1}
$(-1)^{\dim(x)}\mu\left(\hat0,x\right)\ge0$. Hence, $(-1)^jc_j\ge0$ for all
$0\le j<n$. Combine this with the fact \eqref{j} that $c_d<0$.}
\item{Let $0\le j\le d$. The truncation $L^{\ge j}_\Aa$ is CL-shellable
because of being hereditary and \cite[Theorem 10.11]{BW2}, and the
length of any maximal chain in $L^{\ge j}_\Aa\cup\{\hat0\}$ has the
same parity as $j$. Hence, by Hall's theorem \eqref{a},
\cite[Proposition 5.7]{BW1} and Proposition \ref{E}:
$(-1)^{j+1}\wt\chi(L^{\ge j}_\Aa)=(-1)^j c_j\ge0$ for all $0\le j\le
d$. Finally, we know \eqref{k} that $c_d\le0$.}
\end{enumerate}
\end{proof}
\begin{remark}
\item Examples of subspace arrangements covered by Theorem \ref{H} are
hyperplane arrangements, the $k$-equal arrangements $\Aa_{n,k}$ and
$\BB_{n,k}$ for {\em even} $k$ \cite{BSag,BW1}, and several of the orbit arrangements
$\Aa_\la$ shown to be CL-shellable by Kozlov \cite{Ko}.

\item Theorem \ref{H} is not true for general subspace arrangements. For
example, take a planar arrangement $\Aa$ of two intersecting lines and
two points not incident with these lines. In the affine version this
has characteristic polynomial $P_\Aa(t)=t^2-2t-1$, and in the
projective version we get $P^*_\Aa(t)=-t-2$. Note that in these cases
$L_\Aa$ is not mod-$2$-pure but has the other two required properties.

\item If the condition ``mod-$2$-pure'' is strengthened to ``pure'' (all
maximal chains have the same length), then Theorem \ref{H} would
remain valid with ``CL-shellable'' relaxed to ``Cohen-Macaulay''. The
proof is similar, using well-known results about Cohen-Macaulay
complexes. This version of the theorem would however not cover
$\Aa_{n,k}$ and $\BB_{n,k}$ which are not in general pure.
\end{remark}
\end{section}

\let\refitem\DHrefitem

\newcommand\hodli{\mathop{\hbox{\vtop{\setbox0=\hbox{holim}%
\dimen0\wd0\box0
\nointerlineskip\hbox to\dimen0{\rightarrowfill}}}}}

\newcommand\Topos{\symb{{\Cal T}opos}}
\newcommand\Top{\symb{{\Cal T}op}}
\newcommand\Sets{\symb{{\Cal S}ets}}
\newcommand\Sh{\symb{Sh}}
\def\et{\hbox{\tiny \'et}}
\candef{HW}

\begin{section}{The $\Q_\ell$-cohomology of subspace arrangements}

We will now consider the computation of the cohomology of a subspace
arrangement, and in particular its \'etale cohomology. Most of the results that we
will prove are well-known in the case of an arrangement over the real or complex
numbers, and are at least partially to be found in the literature (cf.~\cite{Ya})
over a general field (including positive characteristic). The new contribution
is that we keep track of the action of the Galois group of the base field, which
has important arithmetic significance. As the general ``philosophy of weights''
(cf.~\[D3]) would predict, we can use the same argument to get the mixed
Hodge structure in the complex case, a result which seems to be new (except for
the case of hyperplane configurations which is due to Kim \cite{Ki}). In this
paper we will also be concerned only with results on rational cohomology. In
that case one can use the action of the Galois group (or the rational mixed
Hodge structure) on cohomology to very quickly get to the desired result. 

As we are dealing with varieties over arbitrary fields (our main interest being
the case of finite fields) we are forced to deal with \'etale cohomology of
algebraic varieties instead of classical cohomology, since the latter make sense
only over the real or complex numbers. Its construction is based on the
realisation that to define the usual cohomology one needs access not to the
topological space itself but only the category of sheaves on it. Though neither
the topological space underlying a complex algebraic variety nor its category of
sheaves can be constructed algebraically, a category with properties very
similar to this category of sheaves can be constructed in a purely algebraic
fashion. In general this category is most definitely not equivalent to the
category of sheaves on a topological space, and Grothendieck and his
collaborators \[G4] introduced an axiomatisation under the name of
\Definition{topos} that covers both these new
categories and the category of sheaves on a topological space. The category 
associated to an algebraic variety (or more
generally a scheme) then goes under the name \Definition{\'etale topos}. Our
technical results would be most naturally formulated in terms of toposes and
diagrams of them, but in the interest of concreteness we will confine ourselves
to algebraic varieties (and implicitly their \'etale toposes).
\begin{remark}
It should be noted that in the case of a reasonable topological space the
topological space itself can be recovered from the category of sheaves on it,
hence not only is knowledge of the category of sheaves on a (reasonable)
topological space enough to be able to compute its cohomology, it actually is
equivalent to knowledge about the topological space itself. For the reader
interested in details we can add that ``reasonable'' in this context means that
every \Definition{irreducible} (not the union of two non-empty closed subsets)
closed subset is the closure of a unique point -- a condition almost always
fulfilled in practice.
\end{remark}

The construction of \'etale cohomology is quite involved. The original work of
Grothendieck and his collaborators \cite{D4,G4,G4a,G5} is still the only place
where a detailed treatment of its technical properties can be found. The 
monographs \cite{FK} and \cite{Mi} are easier to approach, but deal mainly
with the case of smooth varieties. We will make a thumbnail sketch of how
the $\ell$-adic \'etale cohomology groups $H^i_{\et}(X,\Q_\ell)$
that we shall use are constructed.

The analogy between the \'etale topos and the category of sheaves on the space
underlying a complex algebraic variety is a very close one, though there are
some definite differences. The most important is that the ``\'etale topology''
is not fine enough to capture the ordinary cohomology with integer coefficients;
one has to use finite coefficients. This is not an artifact of the \'etale topos
but depends on the fact that one wants an algebraically defined
cohomology. Consider for instance the fact that the first cohomology group, with
integer coefficients, of $\C^*:=\C\setminus \{0\}$ is {\Z}. This reflects the
fact that there is a non-trivial covering space of $\C^*$ with structure group
{\Z}, given by $\symb{exp}\co\C \to \C^*$. As the exponential function is
transcendental this makes no algebraic sense, whereas the first cohomology group
with $\Z/n\Z$-coefficients describes covering spaces with structure group
$\Z/n\Z$. In the case of $\C^*$ these are described using $n$'th roots, which
are eminently algebraic functions. One is therefore forced, to begin with, to
work with cohomology with finite coefficients. In that case one obtains a theory
very close to the classical topological one. In fact, a basic theorem
\cite[Exp.~XVI, Thm.~4.2]{G4a} says that for any algebraic variety $X$ over
the complex numbers and any finite group $A$ we have a natural isomorphism of
abelian groups $H^i_{\et}(X,A)\cong H^i_{{\rm cl}}(X,A)$, where the subscripts cl (as in ``classical'') and $\et$ denote respectively the 
ordinary cohomology of the topological space
underlying $X$ and the \'etale cohomology. (This isomorphism is in fact induced
by a map from the topos of sheaves on the topological space of $X$ to the
\'etale topos and hence preserves supplementary structures such as cup
products.) 

Having \'etale cohomology for finite coefficients one then {\em
defines}, for an algebraic variety $X$ over an algebraically closed field,
$H^i_{\et}(X,\Z_\ell)$, $\ell$ a prime, to be the inverse limit
$\ili_nH^i_{\et}(X,\Z/\ell^n\Z)$. Using the result above on equality of \'etale and
classical cohomology for complex varieties $X$ and the universal coefficient
theorem, 
one then gets $H^i_{\et}(X,\Z_\ell)\cong H^i_{{\rm cl}}(X,\Z)\bigotimes \Z_l$. 
If $X$ is defined over a
field $\F$ which is not algebraically closed, then (using a not quite standard
notation) we define $H^i_{\et}(X,\Z_\ell)$ to be the \'etale cohomology of $X$
considered as an algebraic variety over some algebraic closure of the base
field. The fact that $X$ is defined over $\F$ is then reflected in the fact that
we have a natural action of the Galois group of $\F$ on $H^i_{\et}(X,\Z_\ell)$,
of which we will see examples later on. Finally, we put
$H^i_{\et}(X,\Q_\ell):=H^i_{\et}(X,\Z_\ell)\bigotimes_{\Z_\ell}\Q_\ell$, which
then is a finite dimensional vector space over the field $\Q_\ell$. 
We will also normally
dispense with the $\et$-subscript (the comparison theorem guarantees that confusion
should only rarely result.)
\begin{remark}
When the base field is not algebraically closed one may also consider the
\'etale cohomology of $X$ as a variety over the base field. That cohomology will
be an appropriate mixture of the \'etale cohomology of $X$ as an algebraic variety
over an algebraic closure and the Galois cohomology of the base field. As we
will not be interested in it we have chosen to use $H^i_{\et}(X,\Z_\ell)$ to
denote the object which interests us.
\end{remark}

If we consider an arrangement $\Aa$ of subspaces, $V_\Aa$ is by definition their
union. If there are only two of them we would have the Mayer-Vietoris long exact
sequence relating the cohomology of the arrangement, the two linear spaces
covering it and their intersection. In the general case one gets a
Mayer-Vietoris spectral sequence.

The closest analogue of the Mayer-Vietoris long exact sequence would be a
spectral sequence starting with an $E_1$-term. We prefer to start at the
$E_2$-term which, as usual, is more intrinsic. Starting from a covering of an
algebraic variety by closed subvarieties, one may consider the cohomology of
these subvarieties and their intersections. It forms a diagram of abelian groups
over the ordered set of intersections of covering subvarieties. (We will follow
the convention of \[ZZ] in that an ordered set is considered as a category with
morphisms $p \to q$ iff
$p \ge q$, so that a diagram $\{X_p\}$ over the poset has morphisms $X_p \to X_q$
when $p \ge q$). The $E_2$-term will involve the inverse limit and its right
derived functors of this diagram, and we begin by recalling a standard way of
computing such limits.
\begin{lemma}\label{std}
Let \Cal C be a category and $F_.$ a diagram of abelian groups. Then the groups
$\ili_{\Cal C}^*F_.$ are the cohomology groups of the complex $S^*(F)$ whose
$i$'th component is the product
\begin{displaymath}
\prod_{X_0\mapright{f_0}X_1\mapright{f_1} \dots \mapright{f_{i-1}}X_i} F_{X_0}
\end{displaymath}
and whose differential is the alternating sum of the maps obtained by
composing two subsequent morphisms and from the structure map
\pil{f_0^*}{F_{X_0}}{F_{X_1}}.
\pro This is shown in \[JER].
\end{lemma}
\begin{lemma}\label{ss}
Let $\{X_p\}$ be a covering, closed under intersections, of an algebraic variety $X$ by
closed subvarieties. Let $P$ be the poset of these subvarieties ordered
by reverse inclusion. Then there is the Mayer-Vietoris spectral sequence
\begin{equation}
E_2^{i,j} = \ili_{P}{}^j{H^i(X_p,A)} \Rightarrow H^{i+j}(X,A)
\end{equation}
for any finite abelian group $A$.
\pro Let \pil{i_p}{X_p}X be the inclusion map. We may consider the complex
\begin{displaymath}
0 \to A \to \prod_p i_{p*}A \to \prod_{p \ge q}i_{p*}A \to \dots
\end{displaymath}
of sheaves on $X$, where $A$ is considered as the constant sheaf on $X$ and on
the subvarieties $X_p$. To show that this is
an exact sequence it is enough to show that it is exact
on all fibres. For a given point $x \in X$ the fibre at $x$ of this complex is
the (extended) cochain complex with values in A of the (abstract) simplex with
vertices the set of those $p$ for which $x \in X_p$. The simplex being
contractible, this is exact. Now, using that
$H^i(X,i_{p*}A)=H^i(X_p,A)$, as $i_p$ is a closed embedding, and the
spectral sequence of a resolution, we arrive at the Mayer-Vietoris spectral
sequence.
\end{lemma}
We will now apply this result to the rational cohomology of a subspace
arrangement. For this we need to recall some known facts about the action of the
Galois group on \'etale cohomology.  The best control on the Galois action is
obtained when one ignores torsion, so we will look only at \'etale cohomology with
$\Q_\ell$-coefficients (as defined above). It turns out that in positive
characteristic the properties of this cohomology is quite pathological when
$\ell$ is equal to the characteristic of the base field,
so from now on we will assume that {\em
the prime $\ell$ is invertible in the base field}.  
Furthermore, the properties of the Galois
action on cohomology is simplest to formulate when the base field is finite,
so for the moment we will make that assumption and let $q$ denote its
cardinality. Then the Galois group is topologically cyclic (meaning that it has
a dense subgroup generated by one element) with a canonical generator called the
\Definition{Frobenius element}. It is the inverse of the map which raises an
element of an algebraic closure to its $q$'th power.
\begin{remark}
Often it is this map itself rather than its inverse that is called the Frobenius
element, in matters cohomological the present choice is the more suitable
however. The situation is somewhat confusing since the definition of the
Frobenius map in cohomology could appear to give the opposite impression.
However, there is a subtle distinction between the $q$'th power as a generator
of the Galois group of $\F_q$, and the $q$'th power as an algebraic map on, for
instance, affine space. More precisely, both induce actions on the cohomology of
a variety defined over $\F_q$ and these actions are each other's inverses. For
a more thorough discussion of this relation see \cite[pp.~76--81]{D4}.
\end{remark}
We have seen that the action of the Galois group on $\Q_\ell$-cohomology 
of a variety defined over $\F_q$ is
given by a single linear map, usually called the \Definition{Frobenius map}. As
 a first invariant of such a map one may look at the eigenvalues (defined over
 some algebraically closed overfield and counted with the multiplicity in which
 it appears as zeros of the characteristic polynomial). The following definition
 may look very strong.
\begin{definition}
Let $F$ be a linear map of a finite dimensional vector space $V$ over a field of
characteristic zero, and let $q$ be a positive integer. 

\item $F$ is said to be \Definition{pure of weight $n$ (wrt to $q$)} if all of its
eigenvalues are algebraic numbers, all of whose algebraic conjugates have
(complex) absolute value $q^{n/2}$.

\item $F$ is said to be \Definition{mixed} if $V$ has a filtration by $F$-stable
subspaces such that $F$ is pure of some weight on each successive subquotient of
the filtration (where the weight may depend on the subquotient). The set of the
weights of these subquotients will be called the \Definition{weights of $F$}.
\end{definition}
\begin{remark}
\item The condition that all the algebraic conjugates of an algebraic number
have the same absolute value is very strong. For instance, if one bounds
the degree and the absolute value then there are only finitely many such
numbers. This is seen by bounding the coefficients of the minimal polynomial,
coefficients that are also integers.

\item It is implicit in the definition that the set of weights of a mixed linear
operator is independent of the choice of filtration. This is obvious as the set
of weights can be immediately read off from the eigenvalues of $F$.
\end{remark}

A deep result of Deligne (\cite[3.3.8]{D1a}) says that if $X$ is smooth and
proper (over a finite field of cardinality $q$) then its degree $n$
$\Q_\ell$-cohomology is pure of weight $n$, and without any assumptions the
cohomology is mixed.
\begin{example}
\item Affine space is the simplest example, we have
$H^0(\A^n,\Q_\ell)=\Q_\ell$ and the rest of the cohomology groups are equal to
zero. The Frobenius map acts as the identity on this vector space and is hence
pure of weight zero.

\item For the projective line we have $H^1(\P^1,\Q_\ell)=0$, and
$H^0(\P^1,\Q_\ell)=\Q_\ell$, with $F$ acting again as the identity, whereas
$H^2(\P^1,\Q_\ell)$ is more interesting. As a $\Q_\ell$-vector space it is
1-dimensional. The Galois group of the base field acts by the inverse of the
\Definition {cyclotomic character}. The cyclotomic character is the character of
the Galois group for which an element {\gsi} acts by multiplication by the
$\ell$-adic number by which {\gsi} acts on roots of unity of order any power of
$\ell$. In particular, when the base field is finite (of cardinality
$q$) the Frobenius element acts by multiplication by $q$ on
$H^2(\P^1,\Q_\ell)$. Therefore the weight is indeed 2. 

In general we denote by
$\Q_\ell(i)$ the 1-dimensional $\Q_\ell$-vector space on which the Galois group
acts by multiplication by the $i$'th power of the cyclotomic character. Then
from the fact (which can be proved by some of the methods used in the classical
case) that the cohomology of $N$-dimensional projective space is a truncated
polynomial ring and the fact that the Frobenius map preserves multiplication,
one sees that $H^{2i}(\P^N,\Q_\ell)=\Q_\ell(-i)$ when $i \le N$ (and all other
cohomology groups are zero), which confirms that it is indeed pure of weight $2i$.

\item Let us consider the multiplicative group $\mul:=\A^1\setminus \{0\}$. Then
$H^1(\mul,\Q_\ell)=\Q_\ell(-1)$ which is not pure of weight 1 but rather of
weight 2. This does not contradict Deligne's theorem as $\mul$, though smooth, is
not proper. 

More generally, punctured affine space $\A^n\setminus \{0\}$ behaves 
cohomologically as an
odd-dimensional sphere with $H^0(\A^n\setminus \{0\},\Q_\ell)=\Q_\ell$,
$H^{2n-1}(\A^n\setminus \{0\},\Q_\ell)=\Q_\ell(-n)$, and the other
cohomology groups equal to zero.

\item One can also introduce \'etale cohomology with compact support, see below.
Intuitively this
is the reduced cohomology of the one-point compactification --- only that this
doesn't make sense in our setting since the one-point compactification,
even if defined, is not usually an algebraic
variety. For the affine spaces one has that there is only a single non-zero
cohomology group for cohomology with compact support:
$H_c^{2n}(\A^n,\Q_\ell)=\Q_\ell(-n)$.
\end{example}
In all of these examples each weight was associated to only one cohomology
group. In the case of projective space that follows from Deligne's theorem,
in the others it seems to be more of an accident. In any case, they are all
situations to which the following theorem applies. Before going into its
formulation we want to comment on the cohomology of 1-point compactifications, which
occurs prominently in the study of subspace arrangements. 

The 1-point compactification of a complex algebraic variety $X$ need not be an
algebraic variety. If one is interested only in its cohomology a substitute may
be found, the \Definition{cohomology with compact support}. For this one chooses
some realisation of $X$ as an open subset of some proper variety $j\co
X\hookrightarrow{\bar X}$ and then one considers the \'etale cohomology
$H^i(\bar X,j_!A)$, where $j_!A$ is the sheaf on $\bar X$ that is equal to $A$
on $X$ and has fibre 0 at all points of $\bar X$ outside of $X$ (``the extension
by zero''). This turns out to be independent of the choice of $j$ and computes
in the case that the base field is the complex numbers the reduced cohomology of
the 1-point compactification. For these cohomology groups the notation
$H^i_c(X,A)$ is used. (Properly speaking we should also add the subscript 
$\et$, as the cohomology with compact support makes excellent sense also in the
classical case. In the interest of readability we will dispense with that.)
There is now an analogue of the spectral sequence of Lemma \ref{ss} for
cohomology with compact support, the proof is the same (the essential point is
that when one embeds $X$ as an open subset of a proper variety one gets at the
same time a compactification of all the $X_p$ by taking their closures in the
compactification of $X$).
\begin{theorem}\label{degen}
Let $X$ be an algebraic variety that is the union of a family of closed subvarieties
$\{X_p\}$, closed under intersection. Suppose that there is a function {\gph} from
{\N}, the natural numbers, to subsets of the integers such that different numbers are
taken to disjoint sets, and that the degree $i$ cohomology
of each $X_p$ is mixed with weights in $\phi(i)$. Then, with the notations of Lemma
\ref{ss}, the spectral sequence of (loc.~cit.) degenerates at the
$E_2$-term. The same is true if instead cohomology with compact support is
considered.
\pro
Let us first assume that the base field is finite.  If we can prove that
$E_2^{i,j}$ is mixed with weights in $\phi(i)$ we are finished, since then all
the differentials $d^{i,j}_k$ at the $E_k$-term, for $k \ge 2$, will be between
spaces of disjoint weights. However, Lemma \ref{std} presents $E_2^{i,j}$ as a
subquotient of spaces with weights in $\phi(i)$.

For the case of a general base field there are standard techniques for reducing
to the case of a finite base field, for which we refer to for instance
\cite[6.1]{BBD} rather than repeating them here. Very quickly described, one
first uses that base extension from one algebraically closed field to an
algebraically closed overfield does not change cohomology to reduce to the case
where the base field is finitely generated over the prime field. Then there is a
specialisation to a finite field, which again does not change cohomology.
\end{theorem}
\begin{remark}
\item The idea that one could use weights to show that spectral sequences
degenerate is not new. One of its first uses can be found in \[D2], where
it is applied to the study of the cohomology of the complement of a divisor with
normal crossings in a smooth and projective variety. However, there one is
using the mixed Hodge structure on classical cohomology rather than the Galois
action on \'etale cohomology. The arguments of (loc.~cit.) were one of the major
inspirations for our theorem.

\item The theorem applies to the cohomology (including cohomology with compact
support) of an affine subspace arrangement, as there only the cohomology of
affine spaces are involved and we have seen that they fulfill the required
condition. It also applies to projective subspace arrangements, again the
cohomology of projective spaces fulfills the condition. Another case is a
punctured central arrangement, where one considers the arrangement minus a
central point (this is the algebraic analogue of the spherical arrangement
associated to a central arrangement over the reals).
\end{remark}
To apply this result to the various cases of subspace arrangements
we need to compute the cohomology of some diagrams of abelian groups.
Recall from Section 2 the definition of the order complex $\Delta(P)$
of a poset $P$, and of subposets of type $P_{<p}$ and $P_{\le p}$.

\begin{proposition}\label{cohdiag}
Let $P$ be a finite poset and let $A$ be an abelian group.


\item[ii] \label{cohdiag:ii}Let $Q$ be an order ideal (i.e.,
a subset of $P$ such that any element of $P$ less than an
element of $Q$ is also in $Q$), and let $\Cal F_{A,Q}$ be the diagram which is 0
outside of $Q$ and constant with value $A$ on $Q$. Then we have
a natural isomorphism
\begin{displaymath}
\ili{}^j\Cal F_{A,Q} \cong H^j(\Delta(Q),A).
\end{displaymath}

\item[iii] Let $p \in P$ and let $\Cal F_{A,p}$ be the diagram with value $A$ on $p$
and 0 elsewhere. Then we have a natural isomorphism
\begin{displaymath}
\ili{}^j\Cal F_{A,p} \cong \wt H^{j-1}(\Delta(P_{<p}),A).
\end{displaymath}
(For this formula, recall that the reduced cohomology of the empty complex is $A$ in
degree $-1$ and 0 otherwise.)
\pro For part \refitem{ii} we simply use Lemma \ref{std}, which shows that the
higher inverse limits can be computed using a complex which is also the cochain
complex of $\Delta(Q)$ with values in $A$. As for \refitem{iii}, we have a
natural inclusion $\Cal F_{A,p}\hookrightarrow \Cal F_{A, P_{\le p}}$,
whose quotient is $\Cal F_{A, P_{<p}}$. Using the long exact sequence of
higher inverse limits, part \refitem{ii} and the fact that
$\Delta(P_{\le p})$ is contractible, we immediately reach the desired conclusion.
\end{proposition}
\begin{example}
\item Consider the cohomology of an affine arrangement $\Cal A$. The only
non-trivial cohomology group of affine spaces is $H^0(-,A)=A$, so the spectral
sequence degenerates to the isomorphism 
$H^i(V_\Aa,A) \cong H^i(\Delta(L_{\Aa}\setminus\{\hat0\}),A)$.

\item If $\Cal A$ is central we may remove the central point to get an
arrangement of punctured affine spaces (in the real or complex case it is
homotopic to the associated spherical arrangement). Again the condition of Theorem
\ref{degen} is fulfilled. Furthermore, $H^{2i}(-,\Q_\ell)$ is zero for $i>0$ and
$\Q_\ell$ for $i=0$, and $H^{2i-1}(-,\Q_\ell)$ is $\Q_\ell(-i)$ on 
$i$-dimensional elements of the
intersection lattice and zero otherwise. This then is a direct sum
of diagrams of the type considered in the proposition. Hence,
letting $P=L_{\Aa}\setminus \{\hat0\}$ we get
\begin{equation}\label{centrallim}
\begin{array}{lcl}
\ili{}^j H^0(-,\Q_\ell)&=&H^j(\Delta(P),\Q_\ell)\\
\ili{}^j H^{2i-1}(-,\Q_\ell)&=&\Dsum_{\dim p = i}\wt
H^{j-1}(\Delta(P_{<p}),\Z)\bigotimes \Q_\ell(-i).
\end{array}
\end{equation}

\item If $\Cal A$ is a projective arrangement, then $H^{2*+1}(-,A)=0$ and
$H^{2i}(-,A)$ is the constant diagram on the elements of dimension greater than
or equal to $i$. Thus, we may again use the proposition to compute the higher
inverse limits and get the result 
\begin{equation}\label{projlim}
\ili{}^j H^{2i}(-,\Q_\ell)=H^j(\Delta(P^{\ge i}),\Z)\bigotimes \Q_\ell(-i),
\end{equation}
where $P=L_{\Aa}\setminus \{\hat0\}$ and $P^{\ge i} 
:=\set{p \in P}{dim(p) \ge i}$.

\item Once more let $\Aa$ be an affine arrangement, but this time consider
cohomology with compact support. As has been noted, we get a spectral sequence
also in that case, and from the computation of the cohomology with compact
support of affine space we get that $H^{2*+1}_c(-,\Q_\ell)=0$ and that
$H^{2i}_c(-,\Q_\ell)$ is $\Q_\ell(-i)$ on $i$-dimensional elements of the
intersection lattice and zero otherwise. As in the central affine case we get 
\begin{equation}\label{complim}
\ili{}^j H^{2i}_c(-,\Q_\ell)=\Dsum_{\dim p = i}\wt
H^{j-1}(\Delta(P_{<p}),\Z)\bigotimes \Q_\ell(-i).
\end{equation}
\end{example}
Even if one sticks to the case of the base field being the complex
numbers there are advantages to considering diagrams of algebraic varieties. For
algebraic varieties over the complex numbers there is an additional structure on
its cohomology alluded to previously --- its {\it mixed Hodge structure}.

To give the definition of this notion we first recall that a \Definition{Hodge
structure of weight $n$} consists of a finitely generated abelian group $H_\Z$
and a decreasing finite filtration $F^m$ of $H_\Z\bigotimes\C$ by complex
sub-vector spaces such that $\bar F^m$, the complex conjugate of of $F^m$ (the
complex conjugation being induced by that of the second factor in the tensor
product), is a complementary subspace to $F^{n-m+1}$. We then recall
\cite[2.3.1]{D2} that a \Definition{mixed Hodge structure} is a finitely
generated abelian group $H_\Z$ together with one increasing finite filtration
$W_p$ of sub-vector spaces of $H_\Q:=H_\Z\bigotimes_\Z\Q$ and one
decreasing finite filtration $F^m$ by \C-sub-vector spaces of
$H_\C:=H_\Z\bigotimes_\Z\C$, such that for every $i$ the filtration induced by $F^m$ on
$W_i/W_{i-1}\bigotimes\C$ forms a Hodge structure of weight $i$. The
class (with the obvious morphisms) of mixed Hodge structures form an abelian
category. We also use the term \Definition{set of weights} of a mixed Hodge
structure for the set of integers for which $W_i\ne W_{i-1}$. If instead one
looks at only a rational vector space $H_\Q$ without a choice of $H_\Z$ one
speaks about a \Definition{rational Hodge structure}.
\begin{example}
Let $\Z(i)$ be the mixed Hodge structure with $H_\Z=\Z$, $0=W_{2i-1} \subset
W_{2i}=\Q$, and $0=F^{i+1} \subset F^i=\C$\,; and similarly for $\Q(i)$, a
rational Hodge structure. This notation will be used to
describe the mixed Hodge structures relevant to subspace arrangements
after the next theorem.
\end{example}
There is a very strong analogy between mixed Hodge structure and the action of
the Galois group on \'etale cohomology. Parts of the analogy can actually be
proven --- for instance, if one considers the cohomology of a complex algebraic
variety, then the filtration on rational cohomology induced from the Hodge
structure coincides with the weight filtration with respect to the Galois
action.  We illustrate this analogy by giving another proof of the degeneration
of the Mayer-Vietoris spectral sequence when the base field is the complex
numbers.
\begin{theorem}\label{deg}
Let $X$ be a complex algebraic variety that is the union of a family of closed subvarieties
$\{X_p\}$, closed under intersection. Suppose that there is a function {\gph} from
{\N}, the natural numbers, to subsets of the integers such that different numbers are
taken to disjoint sets, and that the degree $i$
cohomology of each $X_p$ 
has weights, with respect to its mixed Hodge structure, in
$\phi(i)$. Then, with the notation of Lemma \ref{ss}, the spectral sequence of
(loc.~cit.) with $\Q$-coefficients degenerates at the $E_2$-term.
\pro
We first need to prove that the spectral sequence is a spectral sequence of
rational mixed Hodge structures. For this we note another way of constructing
it. Namely, we consider the simplicial complex variety $sX_.$ for which $sX_j$
is the disjoint union of the $X_{i_0}$ over the index set $\{i_0 \ge i_1 \ge
\dots \ge i_j\}$ with the obvious structure maps. The spectral sequence
(cf.~\cite[5.3.3.3]{D2a}) applied to the constant sheaf {\Z} of this simplicial
variety converges to the cohomology of $X$ and has an $E_1$-term which is the
standard complex for computing $\ili_P\{H^i(-,\Z)\}$, and hence gives our
spectral sequence from the $E_2$-term on. According to \cite[8.3.5]{D2a} this is
a spectral sequence of mixed Hodge structures which becomes a spectral sequence
of rational mixed Hodge structures when tensored with $\Q$.

If we can prove that $E_2^{i,j}$ is mixed with weights in $\phi(i)$ we are
finished, since then all the differentials $d^{i,j}_k$, for $k \ge 2$, will be
between rational mixed Hodge structures of disjoint weights. However, 
Lemma \ref{std} presents
$E_2^{i,j}$ as a subquotient of spaces with weights in $\phi(i)$.
\end{theorem}
\begin{remark}
We have the following computations of the mixed Hodge structure on the
cohomology of affine space, puctured affine space and projective space,
completely analogous to the action of the Frobenius on \'etale cohomology:
\begin{displaymath}
\begin{array}{lcl}
H^0(\A^n,\Z)&=&\Z(0),\\
H^0(\A^n\setminus\{0\},\Z)&=&\Z(0),\\
H^{2n-1}(\A^n\setminus\{0\},\Z)&=&\Z(n),\\
H^{2i}(\P^n,\Z)&=&\Z(i),\;i \le n.
\end{array}
\end{displaymath}
Hence the theorem may be applied to subspace arrangements.
\end{remark}
Having developed the necessary general tools, we now want to collect our results
as applied to subspace arrangements over finite fields. In that case one extra
refinement is possible which is given in the following lemma.
\begin{lemma}\label{caniso}
Let assumptions be as in Theorem \ref{degen} and assume that the base field is 
finite. Then there is a canonical isomorphism between $H^*(X,\Q_\ell)$
and the $E_2$-term of the spectral sequence. This isomorphism preserves 
the action of the fundamental group of the base field.
\pro We may use the action of the Frobenius map to split up $H^*(X,\Q_\ell)$ as a
sum of generalised eigenspaces under it. Since each such eigenspace occurs in
just one row of the $E_2$-term we get the canonical isomorphism.
\end{lemma}
\begin{remark}
It is not possible to conclude from what we have proven so far that this result
remains true for a general field or has an analogue for mixed Hodge
structures. The reason for this is that it would be possible for the extensions
provided by the spectral sequence to be non-trivial, there are indeed
non-trivial extensions between the Galois representations (resp.~mixed Hodge
structures) involved. It will be proved elsewhere that in the case of subspace
arrangements these possibilities are not realised and in fact the isomorphisms
of the theorem exist for $\Z_\ell$-cohomology (resp.~for cohomology with its mixed
Hodge structure).
\end{remark}
We now collect the various results obtained so far about the cohomology of
unions $V_{\Aa}$ of subspace arrangements over finite fields. To simplify
statements of formulas in this and the following theorem we introduce, 
just as in the classical case, reduced $\ell$-adic cohomology
$\ti H^*(X,\Q_\ell)$. This differs from ordinary cohomology 
for all varieties $X$ only in one dimension, namely
$\ti H^0(X,\Q_\ell)=H^0(X,\Q_\ell)/\Q_\ell$ when $X$ is non-empty, and
$\ti H^{-1}(X,\Q_\ell)=\Q_\ell$ when $X$ is empty.

\newpage
\begin{theorem}\label{coharr}
Let $\Aa$ be a subspace arrangement over a finite field, \pil d{L_\Aa}{\Z}
the dimension function of its intersection semilattice and $P := L_{\Aa}\setminus
\{\hat0\}$. 

\item If $\Cal A$ is an affine arrangement then we have a canonical isomorphism
\begin{displaymath}
H^*(V_\Aa,\Z_\ell)\cong H^*(\Delta(P),\Z_\ell),
\end{displaymath}
which respects the action of the Frobenius map if it is assumed to
act trivially on the right hand side.

\item\label{affcpct} 
If $\Aa$ is an affine arrangement then we have a canonical isomorphism
\begin{displaymath}
H^*_{c}(V_\Aa,\Q_\ell)\cong \bigoplus_{p \in P}
\wt H^{*-2d(p)-1}(\Delta(P_{<p}))\bigotimes\Q_\ell(-d(p)),
\end{displaymath}
which respects the Frobenius action if it is assumed to act trivially on the
cohomology of the order complexes $\Delta(P_{<p})$.

\item If $\Aa$ is a central arrangement then we have a canonical isomorphism
\begin{displaymath}
\ti H^*(V_\Aa\setminus \{0\},\Q_\ell)\cong\bigoplus_{p \in P}
\wt H^{*-2d(p)}(\Delta(P_{<p}))\bigotimes\Q_\ell(-d(p)),
\end{displaymath}
which respects the Frobenius action if it is assumed to act trivially on the
cohomology of the $\Delta(P_{<p})$.

\item\label{projarr} If $\Aa$ is a projective arrangement then we have a
canonical isomorphism
\begin{displaymath}
H^*(V_\Aa,\Q_\ell)\cong
\bigoplus_{0 \le j} H^{*-2j}(\Delta(P^{\ge j}))\bigotimes\Q_\ell(-j),
\end{displaymath}
where $P^{\ge j}:=\set{p \in P}{d(p) \ge j}$, which respects the Frobenius action if
it is assumed to act trivially on the cohomology of the $\Delta(P^{\ge j})$.
\begin{proof}
This follows from the results \ref{ss}, \ref{degen}, \ref{cohdiag} and
\ref{caniso}, together with (\ref{centrallim}), (\ref{projlim}) and
(\ref{complim}).
\end{proof}
\end{theorem}

\medskip
Finally we also collect the consequences of our results for the cohomology of the complements of subspace arrangements over finite fields.
\begin{theorem}\label{cohcomp}
Let $\Aa$ be a subspace arrangement in a space of $n$ dimensions over a
finite field, \pil d{L_\Aa}{\Z} the dimension function of its intersection
semilattice, and $M_\Aa$ the complement of the union $V_\Aa$. Furthermore,
let $P=L_{\Aa}\setminus \{\hat0\}$.

\item If $\Aa$ is affine we have canonical isomorphisms
\begin{displaymath}
H^*_c(M_\Aa,\Q_\ell)\cong \bigoplus_{p \in P}\wt H^{*-2d(p)-2}
(\Delta(P_{<p}))\bigotimes\Q_\ell(-d(p)),
\end{displaymath}
when $* \ne 2n$ and $H^{2n}_c(M_\Aa,\Q_\ell)\cong \Q_\ell(-n)$,
and 
\begin{displaymath}
\wt H^*(M_\Aa,\Q_\ell)\cong \bigoplus_{p \in P}\wt H_{2n-*-2d(p)-2}
(\Delta(P_{<p}))\bigotimes\Q_\ell(d(p)-n).
\end{displaymath} 

\item If $\Aa$ is projective we have canonical isomorphisms
\begin{displaymath}
H^*_c(M_\Aa,\Q_\ell)\cong \bigoplus_{0 \le j \le n}\wt H^{*-2j-1}
(\Delta(P^{\ge j}))\bigotimes\Q_\ell(-j),
\end{displaymath}
when $* \ne 2n$ and $H^{2n}_c(M_\Aa,\Q_\ell)\cong \Q_\ell(-n)$,
and
\begin{displaymath}
\wt H^*(M_\Aa,\Q_\ell)\cong \bigoplus_{0 \le j \le n}\wt H_{2n-*-2j-1}
(\Delta(P^{\ge j}))\bigotimes\Q_\ell(j-n).
\end{displaymath}
All these isomorphisms respect the action of Frobenius if it is assumed
to act trivially on the order complexes occurring in the right hand sides.
\pro One enjoyable property of the cohomology with compact support is
``additivity'' for a closed subvariety and its complement. More precisely, for a
variety $Y$, a closed subvariety $F$ and its complement $U$ we have
(cf.~\cite[Exp.~XVII, 5.1.16.3]{G4a}) a long exact sequence 
\begin{eqnarray*} 
0 \to H^0_c(U,\Q_\ell) \rightarrow &H^0_c(Y,\Q_\ell) \rightarrow& H^0_c(F,\Q_\ell) \rightarrow \\ 
H^1_c(U,\Q_\ell) \rightarrow &H^1_c(Y,\Q_\ell) \rightarrow& H^1_c(F,\Q_\ell) \rightarrow \dots, 
\end{eqnarray*} 
where the maps $H^i_c(Y,\Q_\ell) \to H^i_c(F,\Q_\ell)$ are the restriction
maps. If we apply this to an affine arrangement, using the computation of the
cohomology of affine space as well as that of the union $V_\Aa$, we get that
$H^{2n}_c(M_\Aa,\Q_\ell)=H^{2n}_c(\A^n,\Q_\ell)=\Q_\ell(-n)$ and that
$H^i_c(M_\Aa,\Q_\ell)=H^{i-1}_c(V_\Aa,\Q_\ell)$ for $i \ne 2n$. Now, $M_\Aa$ is
a smooth variety and so we may apply the Poincar\'e duality theorem
\cite[Exp. XVIII,3.2.6.1]{G4a}, which says that the cup-product
$H^i(M_\Aa,\Q_\ell)\bigotimes_{\Q_\ell} H^{2n-i}_c(M_\Aa,\Q_\ell)\to
H^{2n}_c(M_\Aa,\Q_\ell)$ composed with the trace map $H^{2n}_c(M_\Aa,\Q_\ell)\to
\Q_\ell(-n)$ gives a perfect pairing. This gives that $H^i(M_\Aa,\Q_\ell)$ is
canonically isomorphic to $H^{2n-i}_c(M_\Aa,\Q_\ell)^*\bigotimes
\Q_\ell(-n)$. Using this formula, the relation
$H^i_c(M_\Aa,\Q_\ell)=H^{i-1}_c(V_\Aa,\Q_\ell)$ for $i \ne 2n$,
Theorem \ref{affcpct} and the universal coefficient formula applied to the
cohomology of the $\Delta(P_{<p})$, we get the first part of the theorem.

As for the second, we consider again the long exact sequence of cohomology with
compact support, using that for a proper variety it is equal to cohomology
without compact support, so that we can use Theorem \ref{projarr}. Now, it is
clear that the restriction map $H^{2i}(\P^n,\Q_\ell) \to H^{2i}(V_\Aa,\Q_\ell)$
maps $\Q_\ell(-i)$ to $1\otimes \Q_\ell(-i)\subseteq H^0(\Delta(P^{\ge
i}))\otimes \Q_\ell(-i)$.  This is evidently an injection when $i \le m$, where
$m$ is the maximal dimension of subspaces in $\Cal A$, so the long exact
sequence splits up into the desired isomorphisms for cohomology with compact
support. Using duality gives the formula for cohomology without compact support.
\end{theorem}
\begin{remark}
Analogs of the formulas in Theorems \ref{coharr} and \ref{cohcomp} for arrangements
over the real and complex numbers were proved by Goresky and MacPherson \cite{GM},
Ziegler and \v{Z}ivaljevi\'c \cite{ZZ} and others. Some of these formulas in \'etale
cohomology version appear in the paper by Yan \cite{Ya}, however without the
decomposition into eigenspaces under Frobenius.
\end{remark}
\end{section}

\begin{section}{Arrangements over the integers}

In this section we shall be concerned with arrangements specified by
integer forms. Let a $\Z$-{\em arrangement} ({\em affine} resp.~{\em
projective}) mean an arrangement $\Aa=\{K_1,\dots,K_t\}$ where each
subspace is specified by a certain collection of linear forms (general
resp.~homogenous) with integer coefficients. Thus, a $\Z$-arrangement
is really a list of linear forms over $\Z$ partitioned into $t$
groups. With a $\Z$-arrangement $\Aa$ we associate on the one hand the
complex subspace arrangement $\Aa_\C$ (affine or projective, as the
case may be) obtained by interpreting the given $\Z$-forms over $\C$;
and on the other hand the subspace arrangement $\Aa_q$ over the finite
field $\F_q$ obtained from the $\Z$-forms by reduction modulo $p$, for
arbitrary prime powers $q=p^\al$.
\begin{remark}
We could here equally well replace {\Z} with an arbitrary number
ring. Except for trivial notational changes nothing in the arguments to follow
would need to be modified.
\end{remark}
\begin{Lem}\label{A} Let $\Aa$ be a $\Z$-arrangement and $p$ a
prime. Let $\vep$ be the identity map on the set of subspaces of
$\Aa$. Then the following conditions are equivalent:
\begin{enumerate}
\item{$\vep$ extends to a dimension-preserving isomorphism
$L_{\Aa_\C}\cong L_{\Aa_p}$;}
\item{$\vep$ extends to a dimension-preserving isomorphism
$L_{\Aa_\C}\cong L_{\Aa_{p^\al}}$, for all $\al\ge1$;}
\item{$\text{\em{rank}}_\C\{\ell_1,\dots,\ell_g\}=\text{\em{rank}}_{\F_p}\{\ell_1,\dots,\ell_g\}$
for any collection $\ell_1,\dots,\ell_g$ of linear forms from $\Aa$,
containing for each subspace either all of its defining forms or none
of them.}
\end{enumerate}
\end{Lem}

\bpf
The implications (ii)$\Raw$(i)$\Raw$(iii) are immediate. For
(iii)$\Raw$(ii) one checks that the linear algebra in $\F_{p^\al}$ of
the given forms (reduced modulo $p$) takes place in the subfield
$\F_p$.
\epf

We shall call a prime $p$ {\em good} with respect to a
$\Z$-arrangement $\Aa$ if it satisfies the conditions of the lemma,
otherwise {\em bad}. Part (iii) shows that for a given $\Aa$ there is
only a finite number of bad primes (these being the divisors of a
finite collection of determinants in the $\ell_i$'s). In the special
case when $\Aa$ is a hyperplane arrangement condition (iii) can be
expressed by saying that $\Aa$ determines the same matroid over $\C$
and over $\F_p$.

\begin{example} The $k$-equal arrangements defined in Section 2 are
$\Z$-arrangements, and $\Aa_{n,k}$ has no bad primes, while
$\BB_{n,k}$ and $\DD_{n,k}$ have the bad prime 2.
\end {example}
\medskip

Let $\Aa$ be a $d$-dimensional projective $\Z$-arrangement and
$q=p^\al$, where $p$ is a good prime. Let $L_\Aa=L_{\Aa_\C}\cong
L_{\Aa_q}$ and $L^{\ge j}_{\Aa}=\{x\in L_\Aa\mid\dim(x)\ge
j\}\ssm\{\hat0\}$. Define
\begin{equation}\label{m}
\be^{\ge j}_i=\dim_\Q H_i \left(L^{\ge j}_\Aa, \Q\right).
\end{equation}
These order homology Betti numbers of the $j$-truncated intersection
lattices are possibly nontrivial only in the range $0\le i\le d-j\le
d$. We will call the triangular array $(\be^{\ge j}_i)$ the {\em
beta triangle} of $\Aa$.

A formula of Ziegler and \v{Z}ivaljevi\'c \cite[Prop. 2.15]{ZZ} 
\cite[Coroll. 6.7]{WZZ}, which is the complex analog of
Theorem \ref{projarr}, shows that
\begin{equation}\label{o}
\be^\C_i:=\dim_\Q H_i \left(V_{\Aa_\C}, \Q\right)=\sum^d_{j=0} \be^{\ge
j}_{i-2j},
\end{equation}
and from formula \eqref{i} we have that
\begin{equation}\label{p}
Z\left(V_{\Aa_q};t\right)=\prod^d_{j=0}\left(1-q^jt\right)^{\sum(-1)^{i+1}
\be^{\ge j}_i}.
\end{equation}
Thus both the rational Betti numbers of the union of the complex arrangement
$\Aa_\C$ and the zeta function of the discrete arrangement $\Aa_q$ are
governed by the same primitive combinatorial data, namely the beta triangle of
$\Aa$.

\begin{example}
Here is the beta triangle $(\be^{\ge j}_i)$ of
$\Aa_{6,3}$ in the $(i,j)$ Cartesian plane:
\begin{equation*}
\begin{matrix}
20 \\
1 & 26\\
1 & 10 & 10\\
1 & 0 & 0 & 0
\end{matrix}
\end{equation*}
It follows that $\Aa_{6,3}$ has Betti numbers
$\be^\C=(1,\,0,\,1,\,10,\,11,\,26,\,20)$ as a complex variety,
and zeta function
$$
Z\left(V_{\Aa_{6,3}};t\right)=\frac{\left(1-q^2t\right)^{25}}
{\left(1-t\right)\left(1-qt\right)\left(1-q^3t\right)^{20}}.
$$
as a variety over $\F_q$. Furthermore, from Theorem 1.1 we have that
\begin{align*}
P_0(t) &= 1-t\\
P_1(t) &=1\\
P_2(t) &=1-qt\\
P_3(t) &=(1-qt)^{10}\\
P_4(t) &=(1-qt)^{10}\left(1-q^2t\right)\\
P_5(t) &= \left(1-q^2t\right)^{26}\\
P_6(t) &=\left(1-q^3t\right)^{20}
\end{align*}
\end{example}
We will now show that for an important class of $\Z$-arrangements the
Betti numbers $\be^\C_i$ of the complex variety and the zeta function
of the $\F_q$-variety determine each other. This is clearly not true
in general.

The intersection semilattice $L_\Aa$ is said to be {\em
rationally Cohen-Macaulay}  if for all $x<y$ in
$\wh{L}_\Aa=L_\Aa\cup\{\hat1\}$, where $\hat1$ is a new top element,
we have
$$
\wt{H}_i\left(\De(x,y),\Q\right)=0
,\qquad\text{for all\,\,}i<\dim\De(x,y).
$$
This definition is via a theorem of Reisner equivalent to the
Cohen-Macaulayness of the Stanley-Reisner ring of $L_\Aa$. See
Stanley \cite{S1} for more about this concept.

We will consider $\Z$-arrangements $\Aa$ whose semilattice $L_\Aa$ is
both Cohen-Macaul\-ay and hereditary (defined in connection with Theorem
3.8). Then every maximal chain in $L_\Aa$ has the form
$x_0>x_1>\dots>x_d>\hat0$ with $\dim(x_i)=i$ for $0\le i\le
d$. Examples are all hyperplane arrangements, many of the orbit
arrangements $\Aa_\la$ shown to the shellable by Kozlov \cite{Ko}, and
the arrangements corresponding to Cohen-Macaulay simplicial complexes
considered in Bj\"orner and Sarkaria \cite{BSar}. The following
generalises the main result of \cite{BSar}.

\begin{Thm}\label{L}
Let $\Aa$ be a $d$-dimensional projective $\Z$-arrangement such that
$L_\Aa$ is Cohen-Macaulay and hereditary, and let $q$ be a power of a good
prime. Then
$$
Z\left(V_{\Aa_q};t\right)=\prod^d_{j=0}\left(1-q^{d-j}t\right)^{(-1)^{j+1}
\be^\C_{2d-j}-\de_j},
$$
where $\de_j=1$ if $j$ is odd and = 0 otherwise.
\end{Thm}

\bpf
We will use the fact \cite[Theorem III.4.5]{S1} that the truncated
posets $L^{\ge j}_\Aa$ are Cohen-Macaulay for all $0\le j\le d$. Thus,
the beta triangle $(\be^{\ge j}_i)$ has internal zeros, and ones along
the $i=0$ boundary:
\begin{equation*}
\be^{\ge j}_i =
\begin{cases}
1, &\quad\text{if\,\,} i=0, \,0\le j<d\\
0, &\quad\text{if\,\,} 0<i<d-j
\end{cases}.
\end{equation*}
Therefore formula \eqref{o} simplifies as follows for $0\le j\le d$:
\begin{equation*}
\be^\C_{2d-j} =
\begin{cases}
\be^{\ge d-j}_j, &\quad\text{if $j=0$ or $j$ is odd}\\
\be^{\ge d-j}_j+1,&\quad\text{otherwise}
\end{cases}.
\end{equation*}
These two formulas imply
\begin{equation*}
\sum_i(-1)^{i+1} \be^{\ge d-j}_i =
\begin{cases}
\be^\C_{2d-j}-1, &\quad\text{if $j$ is odd}\\
-\be^\C_{2d-j}, &\quad\text{otherwise}
\end{cases},
\end{equation*}
which because of formula \eqref{p} is equivalent to the theorem.
\epf

Note that the rational Betti numbers $\be^\C_{2d-j}$ for $0\le j\le d$
appearing in the theorem are the only essential ones, since the structure
of the beta triangle in the Cohen-Macaulay case shows that for 
$0\le j< d$:

\begin{equation*}
\be^\C_{j} =
\begin{cases}
1, &\quad\text{if $j$ is even}\\
0, &\quad\text{otherwise}
\end{cases}.
\end{equation*}

The preceding proof hinges on the very simple structure of the beta
triangle $(\be^{\ge j}_i)$ given by the almost total vanishing of
Betti numbers in the Cohen-Macaulay case. The beta triangle has
simplified structure also for some other arrangements, including the
$k$-equal arrangements $\Aa_{n,k}$ and $\BB_{n,k}$, as we will now show.

Let us say that an intersection semilattice $L_\Aa$ is {\em
mod-$m$-pure} if the lengths of all maximal chains are congruent mod
$m$.

\begin{Thm}\label{M}
Suppose that $L_\Aa$ is hereditary, mod-$m$-pure and
CL-shellable. 
\newline Then $\wt\be^{\ge j}_i=0$, unless $i+j\equiv d$ (mod $m$).
\end{Thm}

\bpf
Let $0\le j\le d$. As in the proof of Theorem 3.8 we conclude that
$L^{\ge j}_\Aa$ is CL-shellable and mod-$m$-pure. Furthermore, since
$L_\Aa$ is hereditary and $\dim \Aa=d$ there is in $L^{\ge j}_\Aa$ a
maximal chain $x_j>x_{j+1}>\dots>x_d$ with $\dim(x_i)=i$ for all $j\le
i\le d$. Hence, all maximal chains of $L^{\ge j}_\Aa$ have lengths
congruent to $d-j$ (mod $m$), and by \cite[Theorem 5.9]{BW1}
$\wt\be^{\ge j}_i\no=0$ is possible only if $i\equiv d-j$ (mod $m$).
\epf

The intersection lattices of $\Aa_{n,k}$ and $\BB_{n,k}$ satisfy 
these conditions
with $m=k-2$. The nontrivial part here is the CL-shellability, which
was shown in \cite{BW1} and \cite{BSag} respectively.

The material of this section has parallels in the affine case. The
results come out in
essentially the same way, and we will not repeat the arguments.
\end{section}
\vfill\eject

\end{document}